\documentclass[12pt,a4paper,french]{article}

\usepackage{graphics}
\usepackage{graphicx}

\catcode`\é=\active \def é{\' e}

\title{Mean position of a particle submitted to a potential barrier}
\author{D. Mercier, V. R\'egnier \thanks{Laboratoire de Mathématiques
Appliquées et de Calcul Scientifique, Institut des Sciences et Techniques de
Valenciennes, Université de Valenciennes et du Hainaut-Cambrésis, Le Mont
Houy, 59313 VALENCIENNES Cedex 9, FRANCE, courriels :
denis.mercier@univ-valenciennes.fr ; Virginie.Regnier@univ-valenciennes.fr}}

\date{~}

\def\R{\mathop{{\rm I}\kern-.2em{\rm R}}\nolimits}
\def\P{\hbox{I\kern-.2em\hbox{P}}}

\def\N{\hbox{I\kern-.2em\hbox{N}}}

\newcommand{\dfrac}[2]{\displaystyle \frac{#1}{#2}}

\newcommand{\cc}[1]{\cal{#1}}

\def\qqs{\ \forall \ }

\catcode`\à=\active \def à{\` a}
\catcode`\â=\active \def â{\^ a}
\catcode`\é=\active \def é{\' e}
\catcode`\è=\active \def è{\` e}
\catcode`\ê=\active \def ê{\^ e}
\catcode`\ù=\active \def ù{\` u}
\catcode`\ô=\active \def ô{\^ o}
\catcode`\ç=\active \def ç{\c c}
\catcode`\î=\active \def î{\^ \i}

\newcommand{\fin}{\,\rule{1ex}{1ex}\,}

\def\R{\hbox{\it I \hskip-5ptR}}
\newtheorem{theo}{Theorem}

\newtheorem{prop}{Proposition}

\begin{document}
\maketitle

\noindent \bf Abstract \rm \\ A one-dimensional Klein-Gordon problem, which is a physical model for a quantum particle submitted to a potential barrier, is studied numerically : using a variational formulation and a Newmark numerical method, we compute the mean position and standard deviation of the particle as well as their time evolution.      

\vspace{1cm}

\noindent \bf Key words \rm Klein-Gordon equation, Newmark method, mean and standard deviation, kernel smoothing, linear regression.\\
\\
\noindent \bf AMS \rm 35A15, 65M06, 65M12, 81Q05, 81Q10.

\section{Introduction} 
\noindent It has been well-known for a few years now that in quantum mechanics a particle can climb up a step even if it has not enough energy a priori and it will be reflected then with a delay. In classical mechanics it would just try and go back to its position. Everything happens here as if the particle could go through a wall (cf. Fig. 5.1 in \cite{Hai2}) ! This phenomenon is called tunnel effect and has been a subject of interest for physicians and mathematicians : 
the delay has been measured by the physicists A. Haibel and G. Nimtz (cf. \cite{Hai}). It seems that it is approximately equal to the reciprocal of their frequency and that it is independent of the shape of the barrier. The Workshop on Superluminal Velocities gathered some recent contributions in Cologne in 1998 to describe transient phenomena concerning 
the tunnel effect (cf. \cite{Hehl}). In particular J. M. Deutch and F. E. Low were interested in 1993 (see \cite{D-L}) in superluminal effects of low-frequency Gaussian wave packets traversing a rectangular potential barrier. \\
T. Hartman calculated in \cite{Hart} the tunnelling time of Gaussian wave packets for one-di\-men\-sion\-al rectangular potential barriers based on the time-dependent Schr\"{o}dinger equation. An analytic expression is given in
terms of the wave number. It is easy to adapt his calculations to obtain an expression for the delay of 
reflection.\\
The main obstacle for a rigorous formulation comes from the fact that T. Hartman calculates phase shifts for monochromatic waves. Clearly, such waves cannot exhibit physically measurable transient features. Instead signals having a 
carrier frequency have to be studied.\\
That is why F. Ali Mehmeti and V. Régnier used wave packets with a narrow frequency band together with a 
solution formula inspired by J. M. Deutch and F. E. Low in their 2003 paper (\cite{Alreg2} and \cite{Alreg3}). 
For simplicity and to cover the relativistic case we chose the one-dimensional Klein-Gordon equation (as we will do in this paper as well) and formulated our results in terms of the energy flow to ensure the physical measurability of the detected transient phenomena. We gave an analytic expression of the delay which is in accordance with the conjecture of Haibel and Nimtz.  \\ 
More recently (cf. \cite{Alreg4}), we got interested in causality. It is well-known that the solution of the classical linear wave equation with compactly supported initial condition and vanishing initial velocity is also compactly supported in a set depending on time : the support of the solution at time t is causally related to that of the initially given condition. Reed and Simon have shown that for a Klein-Gordon equation with (nonlinear) right-hand side $- \lambda u^3$, causality still holds. We proved in \cite{Alreg4} the same property for a one-dimensional Klein-Gordon problem but with transmission and with a more general repulsive nonlinear right-hand side $F$. We also prove the global existence of a solution using the repulsiveness of $F$. In the particular case $F(u) = - \lambda u^3$, the problem is a physical model for a quantum particle submitted to self-interaction and to a potential step.  \\
\noindent Note that this problem is a transmission problem on a simple network : two semi-infinite branches connected at one point. Recent results on partial differential equations on networks, ramified spaces and more generally on multistructures, can be found in \cite{Lum}. \\
\noindent F. Ali Mehmeti solved explicitly the wave equation in a cross-shaped network in a seminar in 1982 (\cite{fam1}). The result was then generalized by V. Régnier for any integer value of $n$ in \cite{Alreg1}. Analogous coefficients (analogous as our reflection and transmission coefficients) appear in a similar way in the construction of prewavelets in a
star-shaped network in D. Mercier's thesis (\cite{DM}). \\
\\
\noindent In this paper, the geometry is simpler (only two branches) but the waves are dispersive. Our aim is to study the time evolution of the particle especially its mean position and the variance of its position. The idea has been suggested by F. Ali Mehmeti to V. R\'egnier who wants to thank him for doing so. \\
In quantum mechanics neither the position of a particle nor its momentum is known exactly as it would be the case in classical mechanics. The Klein-Gordon equation which we will use in the following is satisfied by a complex-valued function such that the square of the modulus of this function is the probability density of the particle. So the mean position and its variance can give us information on the behaviour of the particle. Since the particle is submitted to a potential step, the potential energy is not constant but only piecewise constant and the problem is the linear transmission following one.

\noindent We call
$\Omega_1 = (- \infty ; 0)$, $\Omega_2 = (0 ; + \infty)$, $I= [0 ; T]$ and denote by
$c$ and $a_k$ ($k=1 ; 2$) three strictly positive constants with $a_2 > a_1$. 
We consider the linear transmission problem $(P)$ of finding $u_k: I \times
\overline{\Omega_{k}} \rightarrow \R$ with the time variable $t \in I$ and the
space variable $x \in \Omega_k$, $k=1 ; 2$ satisfying : 

$$ \left \{ \begin{array} {lllll}
(E_{k}) : & \dfrac {\partial^2 u_{k}}{\partial t^2} (t,x) - c^2
\dfrac{\partial^2 u_{k}}{\partial x^2}  (t,x) + a_k(t,x) u_k(t,x) = 0, \forall (t,x)
\in \R^{+*} \times \Omega_k\\
(T_{0}): & u_{i}(t,0)  =  u_{j} (t,0),\qqs t\in I, ( i,j)\in \{1 ; 2 \} \times 
\{ 1;2\}\\
(T_{1}): & \dfrac{\partial u_1}{\partial x} (t,0^-)  =
\dfrac{\partial u_2}{\partial x} (t,0^+)  , \forall t \in I\\
(IC_{1}): & u_{k}(0,\cdot)  =  f|_{\overline \Omega_k}, k=1 ; 2\\
(IC_{2}): & \dfrac {\partial u_{k}} {\partial t} (0,\cdot) =
g|_{\overline \Omega_k}, k=1 ; 2 
\end{array} \right. $$  

\noindent where $f$ and $g$ are defined from $\R$ to $\R$.\\
\\
\noindent The equations $(E_k)$ are Klein-Gordon equations i.e. dispersive wave equations. The function $f$ represents the initial condition and $g$ contains the initial velocity. This problem is analyzed in \cite{Alreg3}. In particular the meaning of the coefficients is explained~:~the propagation velocity of the waves (phase velocity) is $c$ (in the case $a_k = 0$), and the $a_k$'s are coefficients which characterize the dispersion : the shape of the waves changes when time goes by. Then the velocity of wave propagation is not constant but strongly depends upon the frequency (dual variable in the Fourier transform of the signal). Recall that with classical wave equations ($a_k=0$) the waves are only translated since the propagation velocity is constant (equal to c).\\
When $a_1(t,x)=a_1$ and $a_2(t,x)=a_2$, with $a_2 > a_1 > 0$, the model represents a particle submitted to a potential step with height $(a_2 - a_1)$. \\ 
$(T_0)$ and $(T_1)$ express the absence of energy loss at the central node $x=0$ 
(cf. \cite{fam1}, \cite{fam2}, \cite{fam3} and Theorem 3, Section 3 of \cite{Alreg3}). See Section 5 of the same paper for a physical justification and some more details about the model. \\
\\
\noindent A classical way to treat this problem which is the one we used in \cite{Alreg3} is to use the spectral representation of the spatial operator A involved in the problem and defined as follows :  \\
$D(A)=\{v\in\prod_{k=1}^2H^2(\Omega_k)/ v
\hspace{2mm} \mbox{satisfies} \hspace{2mm} (T_0),(T_1)\}$
\\
$A:v\in D(A)\mapsto Av:= \left(-c^2\dfrac{d^2v_k}{dx^2} + a_k v_k 
\right)_{k\in\{1 ; 2 \}}$ in $H:=~\prod_{k=1}^2 L^2(\Omega_k)$\\ and given 
$u(t)=(u_1(t,\cdot),u_2(t,\cdot))$, the above initial boundary value problem can be rewritten as the following abstract problem denoted by $(AP)$ :

$$(AP) \left \{ \begin{array} {lll}
u \in C^2(\R,H), u(t)\in D(A), \qqs t\in \R^+\\
\dfrac{d^2u}{dt^2}(t)+Au(t)=0 \quad \mbox{in} \quad H\\
u(0)= \Phi, \dfrac{du}{dt}(0)= \Psi
\end{array} \right.$$

\noindent where $\Phi=(f_1, f_2) \in D(A)$ with $f_k : \Omega_k \rightarrow \R$ for $k \in \{ 1 ; 2 \}$ and $\Psi=(g_1, g_2) \in D(\sqrt{A})$ with $g_k : \Omega_k \rightarrow \R$ for $k \in \{ 1 ; 2 \}$. . \\ 
The spectral theory is studied in \cite{fam3}, in Theorem 1.5.1. Here $d_j=1, \qqs j \in \{ 1 \cdots n \}$
and for us the $c_j$'s and $a_j$'s have to be interchanged. The operator $A$ is proved to be self-adjoint and the existence and uniqueness of a solution of the problem in $C^2(\R^+, \prod_{k=1}^n L^2(\Omega_k)) \cap C^1(\R^+, D(\sqrt{A})) \cap C^0(\R^+, D(A))$ is stated in \cite{fam3} (in the context of the analysis of the asymptotic time behaviour of the solution). Recall that $D(\sqrt{A}):=\{v\in\prod_{k=1}^2H^1(\Omega_k)/ v \hspace{2mm} \mbox{satisfies} \hspace{2mm} (T_0)\}$. \\
\\
\noindent The solution in the case $g_1 = g_2= 0$ is recalled in the following theorem (cf. \cite{Alreg3}).

\begin{theo} ~~\\
Assume that $f_k : \Omega_k \rightarrow \R$ with $k \in \{ 1 ; 2 \}$, initial data of Problem (P) are such that $(f_1,f_2) \in D(A)$ . 
\\
Then the restriction to $\Omega_1$ of the unique solution of Problem (P) (such that 
$(u_1,u_2)$ belongs to $C^2(\R^+, \prod_{k=1}^2 
L^2(\Omega_k)) \cap C^1(\R^+, D(\sqrt{A})) 
\cap C^0(\R^+, D(A))$) is given, for $(t,x)$ in $\R^{+*} \times \Omega_1$, by :

$$\begin{array} {llll}
u_1(t,x) =  \frac{1}{2 \pi c^2}  \int_{[\sqrt{a_1};+ \infty)} \cos(\omega t) 
\Im \Big( \frac{2 \omega}{K_1(\omega^2)}   \Big( \int_0^{+\infty} f_1(u)e^{-K_1(\omega^2)(u-x)} 
du \Big) \Big) d\omega \\
\\
 + \frac{1}{2 \pi c^2}  \int_{[\sqrt{a_1};+ \infty)} \cos( \omega t) 
\Im \left( \frac{K_1(\omega^2)-K_2(\omega^2)}{K_1(\omega^2)+K_2(\omega^2)}  \frac{2\omega}{K_1(\omega^2)}  
\left( \int_0^{+\infty} f_1(u)e^{-K_1(\omega^2)(u+x)}du  \right) \right) d\omega\\
\\
 + \frac{1}{2 \pi c^2}  \int_{[\sqrt{a_1};+ \infty)} \cos( \omega t) 
\Im \left( \frac{2K_1(\omega^2)}{K_1(\omega^2)+K_2(\omega^2)}  \frac{2\omega}{K_1(\omega^2)}            
\left( \int_{- \infty}^0 f_2(u)e^{K_2(\omega^2)u-K_1(\omega^2)x}du \right) \right) d\omega 
\end{array}$$

\end{theo}

Recall that, for real $\omega$ and $j \in \{ 1 ; 2 \}$

$$K_j(\omega^2) = \left \{ \begin{array} {ll}
\sqrt{\dfrac{a_j - \omega^2}{c^2}} & \mbox{if} \quad \omega^2 \leq 
a_j \\
i \quad  \sqrt{\dfrac{\omega^2 - a_j}{c^2}} & \mbox{if} \quad \omega^2 \geq a_j
\end{array} \right.$$

\noindent Note that the solution is the sum of three terms : the original term, the reflected one which appears due to the discontinuity in the potential and the last one, called "transmitted term", which contains the information on the transmission of the signal from the second branch to the first one. \\
\\
It is the expansion in generalized eigenfunctions of $\cos(\sqrt{A}t) \Phi$ which is the theoretic solution of $(AP)$ when $\Psi$ vanishes. Otherwise the solution has an additional sine term : it is $\cos(\sqrt{A}t) \Phi + (\sqrt{A})^{-1} \sin(\sqrt{A}t) \Psi$. \\
To study the transient behaviour of the particle we need to compute its mean position $M(t)$ and the variance of its position $V(t)$ defined at the time $t$ respectively by~:
$$\left \{\begin{array}{ll}
M(t) = \dfrac{1}{\int_{\R} |u(t,x)|^2 dx} \left( \int_{\R} x |u(t,x)|^2 dx    \right) \\
\\
V(t) = \dfrac{1}{\int_{\R} |u(t,x)|^2 dx} \left( \int_{\R} \left( x - M(t) \right)^2 |u(t,x)|^2 dx    \right) \\
\end{array} \right.$$      
\noindent Using the spectral exact solution of the problem is too complicated. A numerical approach of the problem is thus used in this paper. \\
In Section 2, the numerical scheme based on Newmark method is explained and stability and error estimates are given using the properties of the exact spectral solution. \\
Section 3 is concerned with the physical interpretation of the results obtained with an initial condition which is sufficiently localized in space and frequency. Heisenberg's principle keeps the position and the momentum of a particle from being known with the same precision : the less uncertain is one of them, the more uncertain is the other one. So only compromises based on Gaussian functions can be chosen and that is our choice. \\
Our main result is that the time series of the mean position of the particle has a linear trend and that of the variance has an exponential one. Except from their sinusoidal variations (linked to the wave behaviour of the particle), the mean position grows linearly with time, at least while the particle has not reached the exterior node (reflections will then take place) and the variance (i.e. the uncertainty of the position of the particle) grows in an exponential way.

\section{Numerical scheme}

\subsection{Variational formulation and Newmark method}

The numerical treatment of the problem, which is developed in the following, is an adaptation of Brézis's work about second-order evolution problems (cf. Chapter 8 of \cite{B}, p. 197-204). 

\noindent The variational formulation of Problem $(P)$ is, for $t \in I = [0;T]$ : 

$$\left \{ \begin{array}{lll}
d^2_t (u(t),v)_H + a(u(t),v) = 0, \qqs v \in V\\
u(0)= f  \qquad (f \in V)\\
u_t(0)= g \qquad (g \in V)
\end{array} \right.$$

\noindent where the scalar product $(\cdot ; \cdot)_H$ is the classical one in the Hilbert space $H=L^2((-L;L))$, $V= H^1((-L;L))$ and the bilinear form $a$ is defined, for $(u,v) \in V \times V$, by
$$a(u,v) = \int_{-L}^L \left( c^2 \partial_x u(x) \cdot \partial_x v(x)+ a(x) \cdot u(x) \cdot v(x) \right) dx$$  

\noindent Now the idea is to construct a sequel of spaces $S_0 \subset S_1 \subset S_2 \cdots S_m \subset V$ and to find a solution $u_m \in S_m$ of Problem $(P_m)$ : 

$$\left \{ \begin{array}{lll}
d^2_t (u_m(t),v)_H + a(u_m(t),v)_H = 0, \qqs v \in S_m\\
u_m(0)= f_m  \qquad (f_m \in S_m)\\
u_{m,t}(0)= g_m \qquad (g_m \in S_m)
\end{array} \right.$$

\noindent Since $f_m$ and $g_m$ are chosen to be approximations of $f$ and $g$ (discretization in the space variable $x$), $u_m$ is an approximation of $u$.  

\noindent Now a discretization in the time variable $t$ is done : $\Delta t = T / n_{max}$ and $t_n = n \Delta t$ with $n = 0 ; 1 ; 2 \cdots n_{max}$. The choices of $n_{max}$ and $m$ depend on the wanted precision of calculus. The numerical scheme  is due to Newmark. If $U^n$ is $u_m(t_n)$ and $V^n = u_{m,t}(t_n)$ then

$$\left \{ \begin{array}{llll}
(U^{n+1} - U^n ; v)_H - (\Delta t) \cdot (V^n ; v)_H + (\Delta t)^2 \left[ a \left( \beta U^{n+1} + \left( \dfrac{1}{2} - \beta \right) U^n ; v \right) \right] = 0\\
\\
(V^{n+1} - V^n ; v)_H + (\Delta t) \cdot a \left( \gamma U^{n+1} + (1 - \gamma) U^n ; v \right) = 0\\
U^0 = f_m  \\
V^0 = g_m 
\end{array} \right.$$

\noindent $\beta$ and $\gamma$ are parameters on which conditions will be put further for the scheme to be stable (cf. the following section.) \\
\noindent Now we choose a basis $\{ \Psi_j^m \}_{j=0}^m$ in the space $S_m$. The Gram matrix of this basis is $G_m = \left( \Psi_j ; \Psi_{j'} \right)_{j,j'}$ and the matrix of the bilinear form $a$ is $A_m = \left( a(\Psi_j ; \Psi_{j'}) \right)_{j,j'}$. Denote by $C^n$, $D^n$, $C_m^0$ and $D_m^0$ the vectors containing the coordinates of $U^n$, $V^n$, $f_m$ and $g_m$ in the basis $\{ \Psi_j^m \}_{j=0}^m$, then the scheme is : 

$$\left \{ \begin{array}{llll}
G_m(C^{n+1} - C^n - (\Delta t) D^n) + (\Delta t)^2 A_m \left[ \beta C^{n+1} + \left( \dfrac{1}{2} - \beta \right) C^n \right] = 0\\
\\
G_m (D^{n+1} - D^n) + (\Delta t) A_m \left( \gamma C^{n+1} + (1 - \gamma) C^n \right) = 0\\
C^0 = C_m^0  \\
D^0 = D_m^0 
\end{array} \right.$$

\noindent Since it holds $u_m(t_n,x)=\sum_{j=0}^m C_j^n \psi_j(x)$, the coefficients $C_j^n$ contain the approximated values of $u$ at $x=x_j$, at the time $t_n=nT/n_{max}$ when the $\psi_j$'s are defined by $\psi_j(x_k)=\delta_{jk}$ for $(j,k) \in \{ 0 ; 1 ; 2 ; \cdots ; m \}^2$ with $\delta_{jk} = 1$ if $j=k$ and 0 otherwise.

\subsection{Stability and error estimates}

The choice of adapted values for $\Delta t$, $\gamma$ and $\beta$ brings stability to the Newmark method.

\begin{prop} ~~\\
Suppose that $f \in D(A^{3/2})$ and $g \in D(A)$, $\beta = 1/4$ and $\gamma = 1/2$. Then there exists a unique solution $u \in C^{3-j}([0;T);D(A^{j/2}))$ for $j = 0 ; 1 ; 2 ; 3$ and the Newmark method is unconditionally stable and of order 2. Denoting by $\Pi_m$ the projection on the vector space $S_m$, the error is 
$$\begin{array}{ll}
|u_m(t_n)-u(t_n)| \leq C \{ |u_m^0 - \Pi_m u^0| + |u_m^1 - \Pi_m u^1| \\
\hspace{3.5cm} + |(I-\Pi_m)u(t_n)| + \int_0^{t_n} \left( |(I - \Pi_m) u_{tt}(s)| + \Delta t |u_t^{(3)}(s)| \right) ds \}
\end{array}$$
\end{prop} 

\noindent The operator $A$ and its domain $D(A)$ are defined in the introduction.\\ 
\noindent \bf Proof. \rm The first part of Theorem 8.6.2 of \cite{B} says that, for regular enough $u$, that is, for $u \in C^2([0,T];V) \cap C^3([0,T];H)$, the method is unconditionally stable (i.e. without any condition on $\Delta t$) and of second order if $\beta = 1/4$ and $\gamma = 1/2$. Now Lemma 1.1.6 of \cite{fam3} states that, if $f \in D(A^{(k+1)/2})$ and $g \in D(A^{k/2})$, then the solution $u$ belongs to $C^{k+1-j}([0;T);D(A^{j/2})$ for $j=0;1;2; \cdots k+1$. Since $V = D(A^{1/2})$ and $H = D(A^0)$, $f \in D(A^{3/2})$ and $g \in D(A)$ are enough conditions. \fin

\section{Numerical results and physical interpretation}

\subsection{Choice of the initial conditions}
\noindent For all the numerical computations, the velocity of the signal $c$ is chosen to be equal to 1 and the potentials $a_1(t,x)$ and $a_2(t,x)$ are constant functions. On the first branch we choose a vanishing potential $a_1=0$. \\
The initial condition $f$ is a Gaussian function multiplied by the square function. This choice is justified by the fact that the initial condition is required to be sufficiently localized in space and frequency i.e. the function and its Fourier transform both need to be localized enough. Now the Fourier transform of a Gaussian function with mean $m$ and standard deviation $\sigma$ is a Gaussian function with mean $m$ and standard deviation $1 / \sigma$. That is why Gaussian functions are good candidates for us. Hardy's note about a remark of Wiener even states that they are the best candidates since the only pair of functions $f$ and $g$ (where $g$ is the Fourier transform of $f$) such that both are $O(|x|^m \exp(-0.5 x^2))$ is given by finite linear combinations of Hermite functions for $f$ and $g$ (see \cite{Har}). \\
The multiplication by the square function then makes the function $f$ vanish at zero. Thus the choice $g = -c f'$ for the initial velocity ($u_t(0)$) makes the signal translate to the right until it bumps into the potential step at $x=0$. Let us recall that the solution of the wave equation with initial conditions $u(t,0) = u_0$ and $u_t(t,0) = v_0$ is 
$$\dfrac{1}{2} \left [ u_0(x+ct) + u_0(x-ct) \right ] + \dfrac{1}{2c} \int_{x-ct}^{x+ct} v_0(y) dy$$ 

\noindent It is the superposition of a left-going term $u^-(t,x)$ and a right-going one $u^+(t,x)$ defined by 
$$\left \{ \begin{array}{ll}
u^+(t,x) = \dfrac{1}{2} u_0(x-ct) - \dfrac{1}{2c} \int_0^{x-ct} v_0(y) dy   \\
\\
u^-(t,x) = \dfrac{1}{2} u_0(x+ct) - \dfrac{1}{2c} \int_0^{x+ct} v_0(y) dy   
\end{array} \right.$$

\noindent So with $v_0 = - c u'_0$ and $u_0$ such that $u_0(0)=0$, $u^-$ vanishes and $u^+(t,x) = u_0(x - ct)$.\\
Physically an adapted initial velocity keeps the signal from "falling" and separating itself into two parts going to opposite directions. 

\noindent For the computations, $f$ is the product of the square function by the Gaussian function with mean $-3$ and standard deviation $1$. A constant factor is added so that $f$ represents a probability density i.e. the integral over $\R$ of $f$ is equal to 1. Then it is 
$$f(x) = \dfrac{1}{10 \sqrt{2 \pi}} \quad x^2 \quad \mbox{exp} \left( - \dfrac{1}{2} (x + 3)^2 \right)$$  
\noindent and $u_t(0,x)$ is   
$$g(x) = \dfrac{1}{10 \sqrt{2 \pi}} \quad x (x^2 + 3 x -2) \quad \mbox{exp} \left( - \dfrac{1}{2} (x + 3)^2 \right)$$  
  
\noindent The length of the branch $L$ is chosen to be equal to $60$ which avoids the reflections at the exterior nodes to take place too fast. They would interfere with what happens at the discontinuity $x=0$, which is our point here. Further since the value of $f$ on $\Omega_2$ is negligible, the probability for the particle to be on the right branch is negligible. So the particle is on the left branch initially (around $- 3.5$ where its probability of presence reaches its maximum) and it meets a potential step. 

\noindent Note that $f$ and $g$ are both rapidly decreasing $C^{\infty}$ functions.

\begin{prop} ~~\\
The functions $f$ and $g$ defined from $\R$ to $\R$ by :
$$\left \{ \begin{array}{ll}
f(x) = \dfrac{1}{10 \sqrt{2 \pi}} \quad x^2 \quad \mbox{exp} \left( - \dfrac{1}{2} (x + 3)^2 \right) \\  
g(x) = \dfrac{1}{10 \sqrt{2 \pi}} \quad x (x^2 + 3 x -2) \quad \mbox{exp} \left( - \dfrac{1}{2} (x + 3)^2 \right)
\end{array} \right.$$  

\noindent satisfy $f \in D(A^{3/2})$ and $g \in D(A)$.
\end{prop}

\noindent \bf Proof. \rm 
Theorem 1.5.2 of \cite{fam3} gives the equivalent statements : \\
(i) $u \in D(A^{(k+1)/2})$ \\
(ii) $u$ satisfies $(K_0), (K_1), ..., (K_k)$, where \\
$u$ satisfies $(K_j) \Longleftrightarrow u \in \prod_{i=1}^2H^{j+1}(\Omega_i)$ and 
$$\left \{ \begin{array} {ll}
A^{j/2}u \quad \mbox{satisfies} \quad (T_0), \quad \mbox{if} \quad j \quad \mbox{is even}\\
A^{(j-1)/2}u \quad \mbox{satisfies} \quad (T_0), \quad \mbox{if} \quad j \quad \mbox{is odd}
\end{array}   \right.$$
\noindent Then $u \in D(A) \Longleftrightarrow u \in \prod_{i=1}^2H^{2}(\Omega_i)$ and $u$ satisfies $(T_0)$ and $(T_1)$. Those conditions of continuity and continuity of the first derivative at $x=0$, are clearly satisfied by the function $g$ which is a $C^{\infty}$ function. \\
Now $u \in D(A^{3/2}) \Longleftrightarrow u \in \prod_{i=1}^2H^{3}(\Omega_i)$, $u$ satisfies $(T_0)$ and $(T_1)$ and $Au$ satisfies $(T_0)$. The first two conditions are satisfied once more. As for the third one, it gives the jump of the second derivative at 0 i.e. $$u''(0^+) - u''(0^-) = \dfrac{a_2}{c^2} u(0)$$
\noindent Since the second derivative of the function $f$ is continuous ($f$ is a $C^{\infty}$ function) and since $f(0)=0$, this transmission condition is satisfied by $f$. \fin

\subsection{Interpretation of the different cases}

\noindent In quantum mechanics, the position of a particle is never known with precision. Only its probability density is known and given by the expression $u(t,\cdot)$ at the time $t$, solution of our problem $(P)$. We calculate the mean $M(t)$ and the standard deviation $\sigma(t)$ of the signal $u(t, \cdot)$ at the time $t$. Let us recall that they are defined as 
$$\left \{\begin{array}{ll}
M(t) = \dfrac{1}{\int_{\R} |u(t,x)|^2 dx} \left( \int_{\R} x |u(t,x)|^2 dx    \right) \\
\\
\sigma (t) = \dfrac{1}{\int_{\R} |u(t,x)|^2 dx} \left( \int_{\R} \left( x - M(t) \right)^2 |u(t,x)|^2 dx    \right) \\
\end{array} \right.$$      

\noindent The variations of $M$ and $\sigma $ with respect to $t$ are studied numerically here in different situations. 

\subsubsection{Necessity of a numerical approach}

\noindent In fact an analytic approach could be envisaged in the absence of potential step since the situation would be less complicated. \\
With a potential step ($a_1=0$ and $a_2 >0$), and in the particular case $g_1=g_2=0$, the solution could be rewritten in terms of Fourier transforms as follows (cf. \cite{Alreg4}) :
 
\begin{theo} ~~\\
Assume that $f_k : \Omega_k \rightarrow \R$ with $k \in \{ 1 ; 2 \}$, initial data of Problem (P), described in the introduction, are such that $(f_1,f_2) \in D(A)$, $f_1$ is compactly supported in $(-\infty ; 0)$ and $f_2 \equiv 0$. 
\\
Then the restriction to $\Omega_1$ of the unique solution of Problem (P) (such that 
$(u_1,u_2)$ belongs to $C^2(\R^+, \prod_{k=1}^2 L^2(\Omega_k)) \cap C^1(\R^+, D(\sqrt{A})) 
\cap C^0(\R^+, D(A))$) is given, for $(t,x)$ in $\R^{+*} \times \Omega_1$, by :

$$u_1(t,x) = \dfrac{1}{\pi} \cc{F}\it^{-1}_{\omega \mapsto x} \left[\cos(\sqrt{c^2 \omega^2}t) \cc{F} \it f_1(\omega) \right]$$

\noindent where the complex square root has been defined such that $\sqrt{c^2 \omega^2 - a_2}$ = 
$$ \left \{ \begin{array} {lll}
\sqrt{c^2 \omega^2 - a_2} \quad \mbox{if} \quad \omega \geq \sqrt{\frac{a_2}{c^2}} \\
i \sqrt{a_2 - c^2 \omega^2}  \quad \mbox{if} \quad - \sqrt{\frac{a_2}{c^2}} \leq \omega \leq \sqrt{\frac{a_2}{c^2}} \\
- \sqrt{c^2 \omega^2 - a_2} \quad \mbox{if} \quad \omega \leq - \sqrt{\frac{a_2}{c^2}}
\end{array} \right.$$ 

\noindent Likewise, for $(t,x)$ in $\R^{+*} \times \Omega_2$, $u_2(t,x) =$  
$$\begin{array} {ll}
 - \dfrac{1}{2 \pi} \int_{- \infty}^{\infty} \cos(\sqrt{c^2 \omega^2}t) 
e^{i \left( \frac{1}{c} \sqrt{c^2 \omega^2 -a_2}  \right) x}
\left( \dfrac{2 \omega c}{\omega c + \sqrt{c^2 \omega^2 -a_2}} \right) \cc{F} \it f_1(\omega) d \omega
\end{array}$$

\end{theo}

\noindent So analysing the time evolution of the mean and of the variance is not simple. If the potential is constant ($a_1=a_2=m^2$), the solution is :
$$u_1(t,x) = \dfrac{1}{\pi} \cc{F}\it^{-1}_{\xi \mapsto x} \left[\cos(\sqrt{m^2+ c^2 \xi^2}t) \cc{F} \it f(\xi) \right]$$
\noindent with the initial condition $f$. Rewriting the cosine as half the sum of two complex conjugated exponentials leads to :
$$u_1(t,x) =  \dfrac{1}{2 \pi} \int_{\R} e^{i(\xi x - \omega(\xi)t)} \cc{F} \it f(\xi) d\xi + \dfrac{1}{2 \pi} \int_{\R} e^{i(\xi x + \omega(\xi)t)} \cc{F} \it f(\xi) d\xi$$
\noindent where $\omega(\xi) = \sqrt{m^2+ c^2 \xi^2}$. The expression $(\xi x - \omega(\xi)t)$ is the phase of the first integral and it is stationary when $x/t=\omega'(\xi)$. Moreover 

\begin{equation}
\omega'(\xi)=\dfrac{c^2\xi}{\sqrt{m^2+ c^2 \xi^2}}
\end{equation}

Thus $\omega$ is an increasing function. If the initial condition is chosen to lie in the frequency band $[\xi_1;\xi_2]$ i.e. if the support of $\cc{F} \it f$ is a subset of $[\xi_1;\xi_2]$, then at a fixed time $t$, the solution $u(.,t)$ lies essentially in $[\omega'(\xi_1) t ; \omega'(\xi_2) t]$. Since the center and the amplitude of this interval both grow linearly with $t$, the mean value and the standard deviation of the particle submitted to a constant potential are due to grow linearly with time. This approach is useful to study the time decay rate of the solution of the Klein-Gordon transmission problem. Marshall, Strauss and Wainger were interested in that problem (cf. \cite{MSW}), F. Ali Mehmeti as well (\cite{fam2} and \cite{fam3}).  \\  
\\
Note that the phase of the second integral can be obtained from the first one, changing $t$ into $-t$. It can be seen as a negative-energy particle solution term propagating backwards in time which describes a positive-energy antiparticle solution term propagating forwards in time. This reinterpretation is due to Feynman and explained in \cite{Ait}. In fact negative-energy solutions to the Klein-Gordon equation may appear due to its construction (cf. Section 6.2.4 in \cite{regth}) and Feynman has given them a physical interpretation. \\
\\
In our situation, the potential is not constant ($a_2-a_1>0$) so the time evolution of the mean position and of its standard deviation is not that clear. This justifies the use of a numerical scheme to study the problem.

\subsubsection{Methodology}

\noindent The numerical scheme, due to Newmark and recalled previously in Section 2.1, has been implemented in a $C^{++}$ program. The conjugate gradient method has been used to solve the system. In each case the graph of the solution $u(t,\cdot)$ is produced for different values of $t$, as well as that of the mean position $M$ and standard deviation $\sigma$ (as functions ot the time $t$). \\
\\
\noindent Two aspects are interesting : \\
\begin{itemize}
\item first the behaviour in the neighbourhood of $x=0$. \\ 
How does the particle interact with the barrier it meets ? Different situations have to be looked at thoroughly : since the initial condition is fixed here, the energy of the particle is fixed and depending on the height of the barrier, whether it will have enough energy to climb it up in a classical way (as a ball would do meeting an inclined plane in classical mechanics) or it has not, but will go through the barrier due to the specific laws of quantum mechanics. The limit case is of course particularly interesting. \\
Let us recall how the energy of the particle is defined. The initial condition $f$ is the probability density of the particle and the Fourier transform of $f$ gives the frequency distribution of the signal. The energy of the particle is proportional to its frequency when it is sufficiently localized in frequency. The Fourier transform of $f$ is $\hat{f}$ defined by :
$$\hat{f}(\omega) = Const \cdot (1 - (\omega + 3)^2) \quad \mbox{exp} \left( - \dfrac{1}{2} (\omega + 3)^2 \right)$$  
where $Const \approx 1.111 \cdot 10^{-3}$. Its shape is that of a mexican hat centered at $-3$. So the particle has enough energy to climb up a barrier with height inferior to $9$ (the square of $3$) and the tunnel effect takes place for higher barriers.
    
\item second, the behaviour for large values of $t$. \\
In fact not that large since we have to stop the experiment just before the wave hits the boundaries to avoid additional reflections. Numerically the initial condition has been chosen far enough from the boundary so that this does not happen too quickly. \\ 
Since the mean value and the standard deviation may both have a wave behaviour (the solution has a wave behaviour), we study the trend of their evolution using Dagnelie's kernel smoothing (cf. Section 15.4.2 of \cite{Dagn}). The idea is to keep only the trend of the phenomena removing their local variations as it is done with time series. The effects of seasonality are removed to be able to see the evolution from one year to the following one. \\
Our conjecture is that, for all the shapes of barriers (infinite, finite, piecewise constant or not) we envisage, both trends are linear. That is checked here using a computation of the correlation coefficient between the mean (respectively the standard deviation) and the time. The equations of the least squares straight lines are also calculated. The notations are : $mean = A t + B$ and $r$ is the correlation coefficient of this regression, and $standard \quad deviation = A_1 t + B_1$ where $r_1$ is the correlation coefficient (cf. Table 1 of Section 3.2.3).

\end{itemize}

\subsubsection{Results and physical interpretation}
\noindent The observation of the different graphs of the solution, its mean value and its standard deviation leads to general remarks about the common behaviour of the solution (independently of the shape of the barrier) : \\
\\
\noindent The standard deviation is always constant till the particle hits the barrier and the mean position grows in a perfect linear way. This is due to the absence of dispersion on the left branch : there is a simple translation of the signal without any reshaping. \\
We observe that both the mean position and the standard deviation have a chaotic behaviour when hitting a high barrier, which is of course not surprising. They are both smoother after the shock. \\
In particular the standard deviation decreases a bit before increasing again and keeps increasing then (see Figure 1). Keeping the analogy with a ball being thrown against a wall, we could say that this small decreasing in the standard deviation can be seen as a reshaping of the particle which is somehow crushed, a bit flattened like a rugby ball when it hits the barrier, exactly like a round ball, which would become imperceptibly oval due to the shock.\\
The analogy is not perfect however : the decreasing of the standard deviation is not larger for particles having more initial energy (as could be expected). It is exactly the same one for a fixed height of the barrier. In classical mechanics a ball thrown with higher velocity (and so more kinetic energy) will get more flattened when hitting the wall. In fact the influent parameter here is the height of the barrier : the standard deviation is larger for higher barriers. What happens is linked to the group velocity of the signal. On the left branch it is $c$, and on the second branch it is given by $(1)$ with $m=a_2$. Thus an increasing height of the barrier $a_2$ implies a smaller velocity of the particle on the left branch. So the signal on the left branch keeps moving to the right (reducing its standard deviation, but independently of the height of the barrier) while the right part of the signal moves slowly and so increases its standard deviation slowly (and the higher the barrier the more slowly it increases). The combination of both phenomena makes the decreasing in the standard deviation  bigger for higher barriers. Everything happens as if a higher barrier corresponded to a "tougher" wall.\\
Another way to explain that the decreasing of the standard deviation is not larger for particles having more initial energy is to say that if the height of the barrier is fixed while the initial energy of the particle increases, then the transmitted part increases. A part of the energy is used by the particle to cross the barrier partially and so it is not flattened as much as if it was completely stopped by the barrier.    \\
Note that in all cases, the standard deviation is affected by the barrier sooner than the mean position. The amplitude of its variations is also larger. In fact the mean keeps on increasing (but a bit more slowly) while the standard deviation decreases. The shock seems to have a small impact on the mean position : a peak during less than $0.05$ units of time while the standard deviation is perturbed for about $0.3$ units of time.\\
\\
To finish with the general results, observed for any value of $a_2$, Table 1 gives the correlation coefficient, denoted by $r$ (respectively $r_1$), between the mean $M$ (respectively the standard deviation $\sigma$) and the time $t$ (cf. Section 3.2.2). As it was conjectured in the introduction, both have a linear behaviour for large values of $t$. It was not that clear a priori for small barriers.

\noindent Let us now consider the different cases : 

\begin{itemize}

\bf \item{Infinite barrier with small height $a_2<9$} : \\
\noindent \rm In that case, the energy of the particle is enough for it to climb up the barrier in the classical way and the smaller the barrier, the easier it is for it to climb up. So the mean position which starts at about $-3.5$ (the particle is on the left branch at the beginning), keeps increasing and takes positive values after some value of $t$ (denoted by $t_0$) which depends on the height of the barrier. The higher the barrier is, the larger $t_0$ becomes. In fact the higher the barrier, the larger the reflected part is (cf. the values of the reflection coefficient represented in Figure 2 of \cite{Alreg3}) and if the reflected part is high, it means that a part of the signal is on the left branch, which makes the mean position smaller. That is why $t_0$ is higher. \\ 
\noindent The height of the barrier influences the behaviour of the particle when hitting the barrier : the peak of the mean position is smoother for a smaller barrier and the increasing of the propagation velocity is not so big for a small barrier since the peak is larger. Note that this apparent acceleration of the particle is in contradiction with the decreasing of the propagation of the velocity on the right branch compared to the left one (cf. $(1)$ once more). In fact this acceleration of the mean position takes place exactly when the standard deviation increase in a brutal way so the mean is less representative of the reality. The uncertainty on the position of the particle increases brutally. \\ 
\\
Also note the wave behaviour of the standard deviation after the shock while the mean position decreases perfectly linearly (except for very small barriers where small oscillations take place). The standard deviation increases more quickly for low barriers than for high ones, reaching at the same time ($15$ units of time) about $6.2$ if $a_2 = 2$, about $2.75$ if $a_2 = 2.5$, about $1.56$ if $a_2 = 4$, about $1.05$ if $a_2 = 5$, about $0.71$ if $a_2 = 6$. Note that the oscillations are more and more numerous (from $5$ if $a_2 = 2$ to $7$ if $a_2 = 6$) but less and less high. The signal is quite flat for $a_2=1$ and $a_2=8$. \\
\\
\noindent If the height of the barrier tends to zero, the situation tends to that of the classical wave equation (no dispersion) so the signal tends to translate without changing its shape. \\

\bf \item{Infinite barrier with large height $a_2>9$} : \\
\rm \noindent In that case, the energy of the particle is not sufficient for it to climb up the barrier in the classical way. The tunnel effect takes place to allow the particle to go through the barrier. A previous energetic study has shown that in fact, the particle lingers on in the neighborhood of $x=0$ for some time, called delay of reflection, before it is reflected (cf. \cite{Alreg3}). We observe that, except for a delay, the energy just after the shock is equal to what it was just before. That delay is visualized by a spatial shift (called Goos-H\"anchen shift) in the two-dimensional case of the double-prism experiment of Haibel and Nimtz (cf. \cite{Hai} and \cite{reg}). \\
After the chaotic behaviour due to the shock with the barrier, the mean position decreases regularly and linearly and the symmetry of the curve means that the reflection is total exactly as we have just explained (the slope of the curve after the shock is opposite to that of the curve before the shock). \\
Furthermore, the more energy the particle has (compared to the barrier it meets), so the lower the barrier, the smaller the delay is (cf. Figures 3 and 5 of \cite{Alreg3}). When $a_2$ increases, the propagation velocity on the left branch decreases so we have conjectured that the particle lingers on for a longer time around zero before being reflected and that it is the reason why the delay is bigger. \\
Observing the shape of the mean position curve for $a_2$ varying between 9 and 100, we detect a peak which is larger and larger as $a_2$ increases. About 0.1 units of time if $a_2=9$, about 0.25 if $a_2=15$ and 0.3 for $a_2=100$. For $a_2 = 9$ the mean position increases with velocity $c$ (i.e. 1 in our numerical computations) till the wave front reaches the potential step. Then it keeps increasing but with higher velocity and decreases with the same velocity after a very small time before decreasing with velocity $c$ again. Now if $a_2$ increases, the mean begins with decreasing (from $-0.7$ to $-0.72$ for $a_2=9.2$, from $-0.79$ to $-1.52$ for $a_2 = 12$ and from $-1$ to $-1.7$ for $a_2 = 100$), before having the same behaviour as for $a_2=9$, except for the maximal value which decreases with increasing $a_2$ (about $-0.55$ if $a_2 = 15$ and $-0.8$ if $a_2 = 5000$). The transmitted part is smaller so the mean is smaller as well.\\
This observation confirms our conjecture : the particle tries to go through the wall and the interactions with the barrier make it linger on for a while, hence the delay in the reflection. \\
\\
As for the standard deviation, it has the same behaviour before the shock as what it was for small barriers but afterwards, it becomes constant again with the same value as for $t=0$. This is due to the total reflection and the absence of dispersion (the potential vanishes) on the left branch.

\bf \item{Other shapes of barrier} : \\
\rm

\end{itemize}

\bf Table 1 (Potential step) \rm \\
\\

\begin{tabular}{|c|c|c|c|c|c|}
\hline
& $a_2 =2$ & $a_2 = 6 $ & $a_2 = 9$  & $a_2 = 15$   & $a_2 = 150$    \\
\hline
$A$ &  $-0.9114$  & $-0.9928$ & $-0.9940$   & $-0.9941$ & $-1$   \\
\hline
$B$  &  $-2.4414$   &$-3.5840$ & $-3.7679$ & $-3.9338$& $-4.2629$ \\
\hline
$r$  &   $-0.9999$  & $-0.9999$ & $-0.9999$ & $-0.9999$ & $-1$ \\ 
\hline 
$A_1$ &  $0.3266$   & $0.0212$ & $0.0013$ & $2.4873 \cdot 10^{-6}$ &  $-2.23 \cdot 10^{-6}$ \\
\hline
$B_1$  & $1.5078$     & $0.6682$ & $0.6533$ & $0.6525$ & $0.653$\\
\hline
$r_1$  & $0.9984$     & $0.9945$ & $0.9905$ & $0.9820$ & $-0.998$ \\ 
\hline
\end{tabular}

\newpage

\begin{figure}[h]
\begin{center} \includegraphics[scale=1]{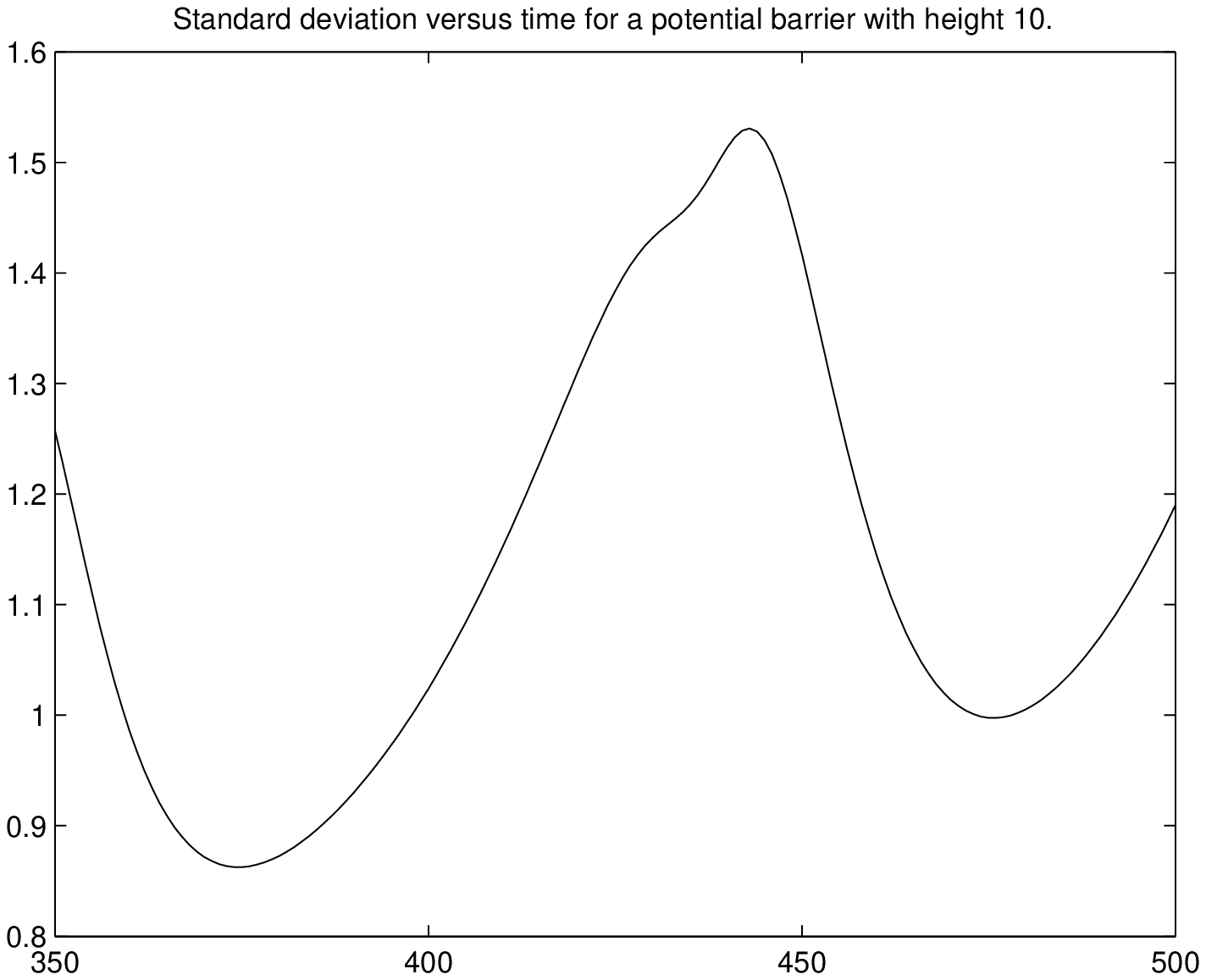}  
\caption{commentaires??\label{figni}}
\end{center}
\end{figure}

\end{document}